# Some Notes on Compact Sets in Soft Metric Spaces

**Sadi Bayramov, Cigdem Gunduz (Aras), Murat I. Yazar**

Department of Mathematics, Kafkas University, Kars, 36100-Turkey
Department of Mathematics, Kocaeli University, Kocaeli, 41380-Turkey
Department of Mathematics, Kafkas University, Kars, 36100-Turkey

baysadi@gmail.com, carasgunduz@gmail.com, miy248@yahoo.com

**Abstract**

The first aim of this study is to define soft sequential compact metric spaces and to investigate some important theorems on soft sequential compact metric space. Second is to introduce $\tilde{\varepsilon}-$ net and totally bounded soft metric space and study properties of this space. Third is to define Lebesque number and soft uniformly continuous mapping and investigate some theorems in detail.

## 1. Introduction

Many practical problems in economics, engineering, environment, social science, medical science etc. cannot be dealt with by classical methods, because classical methods have inherent difficulties. The reason for these difficulties may be due to the inadequacy of the theories of parameterization tools. Molodtsov [12] initiated the concept of soft set theory as a new mathematical tool for dealing with uncertainties. Maji et al. [10] gave some operations on soft sets. Maji et al. [9] introduced some results on an application of fuzzy-soft-sets in decision making problem.
M. Shabir and M. Naz [13] presented soft topological spaces and they investigated some properties of soft topological spaces. Later, many researches about soft topological spaces were studied in [7,8,11,13,14,15,16]. In these studies, the concept of soft point is expressed by different approaches. In the study we use the concept of soft point which was given in [4,14].
Recently in [4,5] S. Das and S. K. Samanta introduced a different notion of soft metric space by using a different concept of soft point and investigated some basic properties of these spaces.

In the present paper, we first give, as preliminaries, some well-known results in soft set theory. We define soft sequential compact metric spaces and investigate some important theorems on soft sequential compact metric space. We also introduce $\tilde{\varepsilon}-$ net and totally bounded soft metric space and study properties of this space. Sequential compactness can be expressed by totally boundedness. Finally, we define Lebesque number and soft uniformly continuous mapping and investigate some theorems in detail.

## 2. Preliminaries

Throughout this paper, $X$ refers to an initial universe, $E$ is the set of all parameters for $X$.

**Definition 2.1. ([12])** A pair $(F,A)$ is called a soft set over $U$, where $F$ is a mapping given by $F: A \to P(U)$.

In other words, the soft set is a parameterized family of subsets of the set $U$. For $\varepsilon \in A$, $F(\varepsilon)$ may be considered as the set of $\varepsilon$-elements of the soft set $(F,A)$, or as the set of $\varepsilon$-approximate elements of the soft set.

According to this manner, a soft set $(F,A)$ is given as consisting of collection of approximations:
$$(F,A) = \{F(\varepsilon) : \varepsilon \in A\}.$$

**Definition 2.2. ([10])** For two soft sets $(F,A)$ and $(G,B)$ over $U$, $(F,A)$ is called a soft subset of $(G,B)$ if

$(i)\ A \subset B$, and

$(ii)\ \forall \varepsilon \in A, F(\varepsilon)$ and $G(\varepsilon)$ are identical approximations.

This relationship is denoted by $(F,A) \tilde{\subset} (G,B)$. Similarly, $(F,A)$ is said to be a soft superset of $(G,B)$, if $(G,B)$ is a soft subset of $(F,A)$. This relationship is denoted by $(F,A) \tilde{\supset} (G,B)$.

**Definition 2.3. ([10])** Two soft sets $(F,A)$ and $(G,B)$ over $U$ are said to be soft equal if $(F,A)$ is a soft subset of $(G,B)$ and $(G,B)$ is a soft subset of $(F,A)$.

**Definition 2.4. ([10])** The intersection of two soft sets $(F,A)$ and $(G,B)$ over $U$ is the soft set $(H,C)$, where $C = A \cap B$ and $\forall \varepsilon \in C$, $H(\varepsilon) = F(\varepsilon) \cap G(\varepsilon)$. This is denoted by $(F,A) \tilde{\cap} (G,B) = (H,C)$.

**Definition 2.5. ([10])** If $(F,A)$ and $(G,B)$ are two soft sets, then $(F,A)$ **AND** $(G,B)$ is denoted $(F,A) \wedge (G,B)$. $(F,A) \wedge (G,B)$ is defined as $(H, A \times B)$, where $H(\alpha, \beta) = F(\alpha) \cap G(\beta), \forall (\alpha, \beta) \in A \times B$.

**Definition 2.6. ([10])** The union of two soft sets $(F,A)$ and $(G,B)$ over $U$ is the soft set, where $C = A \cup B$ and $\forall \varepsilon \in C$,
$$H(\varepsilon) = \begin{cases} F(\varepsilon), & \text{if } \varepsilon \in A - B \\ G(\varepsilon), & \text{if } \varepsilon \in B - A \\ F(\varepsilon) \cup G(\varepsilon), & \varepsilon \in A \cap B \end{cases}.$$

This relationship is denoted by $(F,A) \tilde{\cup} (G,B) = (H,C)$.

**Definition 2.7. ([11])** The complement of a soft set $(F,A)$ is denoted by $(F,A)^c$ and is defined by $(F,A)^c = (F^c, A)$, where $F^c : A \to P(X)$ is a mapping given by $F^c(\alpha) = X - F(\alpha), \forall \alpha \in A$.

**Definition 2.8. ([3])** Let $\mathbb{R}$ be the set of real numbers and $B(\mathbb{R})$ be the collection of all non-empty bounded subsets of $\mathbb{R}$ and $E$ taken as a set of parameters. Then a mapping $F : E \to B(\mathbb{R})$ is called a soft real set. It is denoted by $(F, E)$. If specifically $(F, E)$ is a singleton soft set, then identifying $(F, E)$ with the corresponding soft element, it will be called a soft real number and denoted $\tilde{r}, \tilde{s}, \tilde{t}$ etc.
$\bar{0}, \bar{1}$ are the soft real numbers where $\bar{0}(e) = 0, \bar{1}(e) = 1$ for all $e \in E$, respectively.

**Definition 2.9. ([3])** For two soft real numbers
(i) $\quad \tilde{r} \leq \tilde{s}$ if $\tilde{r}(e) \leq \tilde{s}(e)$, for all $e \in E$;
(ii) $\quad \tilde{r} \geq \tilde{s}$ if $\tilde{r}(e) \geq \tilde{s}(e)$, for all $e \in E$;
(iii) $\quad \tilde{r} < \tilde{s}$ if $\tilde{r}(e) < \tilde{s}(e)$, for all $e \in E$;
(iv) $\quad \tilde{r} > \tilde{s}$ if $\tilde{r}(e) > \tilde{s}(e)$, for all $e \in E$.

**Definition 2.10. ([4,14])** A soft set $(P, E)$ over $X$ is said to be a soft point if there is exactly one $e \in E$, such that $P(e) = \{x\}$ for some $x \in X$ and $P(e') = \varnothing, \forall e' \in E / \{e\}$. It will be denoted by $\tilde{x}_e$.

**Definition 2.11. ([4,14])** Two soft point $\tilde{x}_e, \tilde{y}_{e'}$ are said to be equal if $e = e'$ and $P(e) = P(e')$ i.e., $x = y$. Thus $\tilde{x}_e \neq \tilde{y}_{e'} \Leftrightarrow x \neq y$ or $e \neq e'$.

**Proposition 2.12. ([14])** The union of any collection of soft points can be considered as a soft set and every soft set can be expressed as union of all soft points belonging to it; $(F, E) = \bigcup_{\tilde{x}_e \in (F, E)} \tilde{x}_e$.

**Definition 2.13. ([13])** Let $\tau$ be the collection of soft sets over $X$, then $\tau$ is said to be a soft topology on $X$ if
(1) $\Phi, X$ belong to $\tau$
(2) the union of any number of soft sets in $\tau$ belongs to $\tau$
(3) the intersection of any two soft sets in $\tau$ belongs to $\tau$.

The triplet $(X, \tau, E)$ is called a soft topological space over $X$.

**Definition 2.14. ([8])** Let $(X, \tau, E)$ be a soft topological space over $X$. Then soft interior of $(F, E)$, denoted by $(F, E)^\circ$, is defined as the union of all soft open sets contained in $(F, E)$.

**Definition 2.15. ([8])** Let $(X, \tau, E)$ be a soft topological space over $X$. Then soft closure of $(F, E)$, denoted by $\overline{(F, E)}$, is defined as the intersection of all soft closed super sets of $(F, E)$.

**Definition 2.16. ([7])** Let $(X, \tau, E)$ and $(Y, \tau', E)$ be two soft topological spaces, $f : (X, \tau, E) \to (Y, \tau', E)$ be a mapping. For each soft neighborhood $(H, E)$ of $(f(x)_e, E)$, if there exists a soft neighborhood $(F, E)$ of $(x_e, E)$ such that $f((F, E)) \subset (H, E)$, then $f$ is said to be soft continuous mapping at $(x_e, E)$.

If $f$ is soft continuous mapping for all $(x_e, E)$, then $f$ is called soft continuous mapping.

Let $\tilde{X}$ be the absolute soft set i.e., $F(e) = X, \forall e \in E$, where $(F, E) = \tilde{X}$ and $SP(\tilde{X})$ be the collection of all soft points of $\tilde{X}$ and $\mathbb{R}(E)^*$ denote the set of all non-negative soft real numbers.

**Definition 2.17.** ([4]) A mapping $\tilde{d} : SP(\tilde{X}) \times SP(\tilde{X}) \to \mathbb{R}(E)^*$, is said to be a soft metric on the soft set $\tilde{X}$ if $d$ satisfies the following conditions:

(M1) $\tilde{d}(\tilde{x}_{e_1}, \tilde{y}_{e_2}) \tilde{\geq} \overline{0}$ for all $\tilde{x}_{e_1}, \tilde{y}_{e_2} \tilde{\in} \tilde{X}$,

(M2) $\tilde{d}(\tilde{x}_{e_1}, \tilde{y}_{e_2}) = \overline{0}$ if and only if $\tilde{x}_{e_1} = \tilde{y}_{e_2}$,

(M3) $\tilde{d}(\tilde{x}_{e_1}, \tilde{y}_{e_2}) = \tilde{d}(\tilde{y}_{e_2}, \tilde{x}_{e_1})$ for all $\tilde{x}_{e_1}, \tilde{y}_{e_2} \tilde{\in} \tilde{X}$,

(M4) For all $\tilde{x}_{e_1}, \tilde{y}_{e_2}, \tilde{z}_{e_3} \tilde{\in} \tilde{X}$, $\tilde{d}(\tilde{x}_{e_1}, \tilde{z}_{e_3}) \tilde{\leq} \tilde{d}(\tilde{x}_{e_1}, \tilde{y}_{e_2}) + \tilde{d}(\tilde{y}_{e_2}, \tilde{z}_{e_3})$.

The soft set $\tilde{X}$ with a soft metric $\tilde{d}$ on $\tilde{X}$ is called a soft metric space and denoted by $(\tilde{X}, \tilde{d}, E)$.

**Definition 2.18.** [5] Let $(F, E)(\neq) \in S(\tilde{X})$, then the collection of all soft elements of $(F, E)$ will be denoted by $SE(F, E)$. For a collection B of soft elements of $\tilde{X}$, the soft set generated by B is denoted by $SS(B)$.

**Definition 2.19.** [5] A mapping $\tilde{d} : SE(\tilde{X}) \times SE(\tilde{X}) \to \mathbb{R}(E)^*$, is said to be a soft metric on the soft set $\tilde{X}$ if $d$ satisfies the following conditions:

(M1) $\tilde{d}(\tilde{x}, \tilde{y}) \tilde{\geq} \overline{0}$ for all $\tilde{x}, \tilde{y} \tilde{\in} \tilde{X}$,

(M2) $\tilde{d}(\tilde{x}, \tilde{y}) = \overline{0}$ if and only if $\tilde{x} = \tilde{y}$,

(M3) $\tilde{d}(\tilde{x}, \tilde{y}) = \tilde{d}(\tilde{y}, \tilde{x})$ for all $\tilde{x}, \tilde{y} \tilde{\in} \tilde{X}$,

(M4) For all $\tilde{x}, \tilde{y}, \tilde{z} \tilde{\in} \tilde{X}$, $\tilde{d}(\tilde{x}, \tilde{z}) \tilde{\leq} \tilde{d}(\tilde{x}, \tilde{y}) + \tilde{d}(\tilde{y}, \tilde{z})$.

The soft set $\tilde{X}$ with a soft metric $\tilde{d}$ on $\tilde{X}$ is called a soft metric space and denoted by $(\tilde{X}, \tilde{d}, E)$.

**Note 2.20.** The metric defined in the definition 2.18 is different from the metric given in the definition 2.19 interms of the soft point concept. In this study we use the metric that is given in the definition 2.18.

**Definition 2.21.** ([4]) Let $(\tilde{X}, \tilde{d}, E)$ be a soft metric space and $\tilde{\varepsilon}$ be a non-negative soft real number. $B(\tilde{x}_e, \tilde{\varepsilon}) = \{\tilde{y}_{e'} \tilde{\in} \tilde{X} : \tilde{d}(\tilde{x}_e, \tilde{y}_{e'}) \tilde{<} \tilde{\varepsilon}\} \subset SP(\tilde{X})$ is called the soft open ball with center $\tilde{x}_e$ and radius $\tilde{\varepsilon}$ and $B[\tilde{x}_e, \tilde{\varepsilon}] = \{\tilde{x}_e \tilde{\in} \tilde{X}; \tilde{d}(\tilde{x}_e, \tilde{y}_{e'}) \tilde{\leq} \tilde{\varepsilon}\} \subset SP(\tilde{X})$ is called the soft closed ball with center $\tilde{x}_e$ and radius $\tilde{\varepsilon}$.

**Definition 2.22.** ([4]) Let $\{\tilde{x}_{\lambda,n}\}_n$ be a sequence of soft points in a soft metric space $(\tilde{X}, \tilde{d}, E)$. The sequence $\{\tilde{x}_{\lambda,n}\}_n$ is said to be convergent in $(\tilde{X}, \tilde{d}, E)$ if there is a soft point $\tilde{y}_\mu \tilde{\in} \tilde{X}$ such that $d(\tilde{x}_{\lambda,n}, \tilde{y}_\mu) \to \bar{0}$ as $n \to \infty$.

This means for every $\tilde{\varepsilon} \tilde{>} \bar{0}$, chosen arbitrarily, $\exists$ a natural number $N = N(\tilde{\varepsilon})$, such that $\bar{0} \tilde{\leq} d(\tilde{x}_{\lambda,n}, \tilde{y}_\mu) \tilde{<} \tilde{\varepsilon}$, whenever $n > N$.

**Theorem 2.23.** ([4]) Limit of a sequence in a soft metric space, if exist is unique.

**Definition 2.24.** ([4]) (Cauchy Sequence). A sequence $\{\tilde{x}_{\lambda,n}\}_n$ of soft points in $(\tilde{X}, \tilde{d}, E)$ is considered as a Cauchy sequence in $\tilde{X}$ if corresponding to every $\tilde{\varepsilon} \tilde{>} \bar{0}, \exists m \in N$ such that $d(\tilde{x}_{\lambda,i}, \tilde{x}_{\lambda,j}) \tilde{\leq} \tilde{\varepsilon}, \forall i, j \geq m$, i.e., $d(\tilde{x}_{\lambda,i}, \tilde{x}_{\lambda,j}) \to \bar{0}$ as $i, j \to \infty$.

**Definition 2.25.** ([4]) (Complete Metric Space). A soft metric space $(\tilde{X}, \tilde{d}, E)$ is called complete if every Cauchy Sequence in $\tilde{X}$ converges to some point of $\tilde{X}$. The soft metric space $(\tilde{X}, \tilde{d}, E)$ is called incomplete if it is not complete.

**Definition 2.26.** ([16]) Let $(X, \tau, E)$ be a soft topological space and $W \subseteq X$
(1) A family $C = \{(F_A)_i : i \in J\}$ of open soft sets is called an open cover of X, if it satisfies $\tilde{\bigcup}_{i \in J}(F_A)_i = \tilde{E}$, for each $e \in E$. A finite subfamily of a soft open cover $\{(F_A)_i : i \in J\}$ of $W$ is called a finite subcover of $\{(F_A)_i : i \in J\}$.
(2) $W$ is called soft compact if every soft open cover of W has a finite subcover.

## 3. Compact sets on soft metric spaces

In this section, we study some important properties of soft metric spaces.

**Definition 3.1.** Let $(\tilde{X}, \tilde{d}, E)$ be a soft metric space. $(\tilde{X}, \tilde{d}, E)$ is called soft sequential compact metric space if every soft sequence has a soft subsequence that converges in $\tilde{X}$.

In the following, we give relations between sequential compactness of parameter metric spaces and soft sequential compactness of soft metric spaces.

**Proposition 3.2.** If $(\tilde{X}, \tilde{d}, E)$ is a soft sequential compact, then $(X, d_e)$ is sequential compact, for each $e \in E$.

**Proof.** Let $\left(\tilde{X},\tilde{d},E\right)$ be a soft sequential compact and $\{x^n\}_n$ be any sequence in $(X,d_e)$, for each $e \in E$. By using the sequence $\{x^n\}_n$, soft sequence $\{x_e^n\}_n$ is obtained in $\left(\tilde{X},\tilde{d},E\right)$. Then there exists a soft subsequence as $\{x_e^{n_k}\}_n$. Hence subsequence $\{x^{n_k}\}$ converges in $(X,d_e)$.

The converse of the Proposition 3.2 above may not be true in general. This is shown by the following example.

**Example 3.3.** Let $E = \mathbb{N}$, $X = [0,1]$ and consider the soft metric $\tilde{d}(x_e, y_{e'}) = |e - e'| + |x - y|$, where $\mathbb{N}$ is natural number set. It is clear that $(X, d_e)$ is sequential compact, for each $e \in E$. However, soft sequence $\left\{\left(\frac{1}{2^n}\right)_n\right\}_n$ does not have a convergent soft subsequence in $\left(\tilde{X},\tilde{d},E\right)$.

**Proposition 3.4.** Let $\left(\tilde{X},\tilde{d},E\right)$ be a soft metric space. $\left(\tilde{X},\tilde{d},E\right)$ is soft sequential compact if and only if every infinite soft set $(F,E)$ has a soft cluster point.

**Proof.** $\Rightarrow$ Let $(F,E)$ be infinite soft set and $\{x_{e_n}^n\}$ be a soft sequence in $(F,E)$. Then $\{x_{e_n}^n\}$ has a convergent soft subsequence as $\{x_{e_{n_k}}^{n_k}\}$. Assume that $\{x_{e_{n_k}}^{n_k}\}$ converges to a soft point $z_{e_0}$. Since $\Phi \neq B\left(z_{e_0}, \tilde{\varepsilon}\right) \cap \{x_{e_{n_k}}^{n_k}\} \subset B\left(z_{e_0}, \tilde{\varepsilon}\right) \cap (F,E)$, $z_{e_0}$ is a soft cluster point in $(F,E)$.

$\Leftarrow$: Let $\{x_{e_n}^n\}_n$ be an arbitrary soft sequence. If the sequence is finite, then fixed convergent subsequence of soft sequence can be found. Let $\{x_{e_n}^n\}_n = (F,E)$, then we say that $\{x_{e_n}^n\}_n$ is a soft set. According to condition, this soft set has a soft cluster point as $z_{e_0}$. Since
$$B\left(z_{e_0}, \tilde{k}\right) \cap (F,E) \neq \Phi,$$
we choose $x_{e_{n_k}}^{n_k} \in B\left(z_{e_0}, \tilde{k}\right) \cap (F,E)$, for each $\tilde{k} > \tilde{0}$. Thus the soft subsequence $\{x_{e_{n_k}}^{n_k}\}$ converges to $z_{e_0}$.

**Definition 3.5.** Let $\left(\tilde{X},\tilde{d},E\right)$ be a soft metric space, $S$ be a soft set of soft points. If the condition $\tilde{X} \subset \bigcup_{x_e \in S} B\left(x_e, \tilde{\varepsilon}\right)$ is satisfied, then $S$ is said to be a soft $\tilde{\varepsilon}$-net in $\left(\tilde{X},\tilde{d},E\right)$.

**Definition 3.6.** Let $\left(\tilde{X},\tilde{d},E\right)$ be a soft metric space. If there exists a finite $\tilde{\varepsilon}$-net of $\left(\tilde{X},\tilde{d},E\right)$, for each $\tilde{\varepsilon} > \tilde{0}$, $\left(\tilde{X},\tilde{d},E\right)$ is said to be totally bounded.

**Lemma 3.7.** Let $\left(\tilde{X},\tilde{d},E\right)$ be a totally bounded soft metric space, $A$ be an infinite soft set. Then there is an infinite soft set $B \subset A$ such that $\tilde{d}(B) < \tilde{\varepsilon}$, for each $\tilde{\varepsilon} > \tilde{0}$.

**Proof.** Let $\tilde{\varepsilon} > \tilde{0}$ be arbitrary. Since $\left(\tilde{X},\tilde{d},E\right)$ is totally bounded, there exists a finite $\frac{\tilde{\varepsilon}}{3}$ – net as $H = \left\{x_{e_1}^1, x_{e_2}^2, ..., x_{e_{n1}}^n\right\}$. Then since $\tilde{X} = \bigcup_{i=1}^{n} B\left(x_{e_i}^i, \tilde{\varepsilon}\right)$, $\tilde{A} = \bigcup_{i=1}^{n} B\left(x_{e_i}^i, \tilde{\varepsilon}\right) \cap A$ is obtained. For $1 \leq i \leq n$, at least one of the soft sets $\tilde{A} = \bigcup_{i=1}^{n} B\left(x_{e_i}^i, \tilde{\varepsilon}\right) \cap A$ must have infinite elements. If we denote this set by $B$, it is clear that $\tilde{d}(B) < \tilde{\varepsilon}$.

**Theorem 3.8.** Let $\left(\tilde{X},\tilde{d},E\right)$ be a soft metric space. $\left(\tilde{X},\tilde{d},E\right)$ is totally bounded if and only if every soft sequence has a Cauchy soft subsequence in $X$.

**Proof.** $\Rightarrow$ Let $\left(\tilde{X},\tilde{d},E\right)$ be totally bounded and $\left\{x_{e_n}^n\right\}_n$ be any soft sequence. If the soft sequence $A = \left\{x_{e_n}^n\right\}_n$ is finite, the proof is completed. Assume that $A = \left\{x_{e_n}^n\right\}_n$ is infinite. From the Lemma 3.7, there exists an infinite soft set $B_1 \subset A$ such that $\tilde{d}(B_1) < \tilde{1}$. We choose number $n_1$ such that $x_{e_{n_1}}^{n_1} \in B_1$. If we apply Lemma 3.7 to $B_1$, we obtain an infinite $B_2 \subset B_1$ such that $\tilde{d}(B_1) < \frac{\tilde{1}}{2}$. Here we take the number $n_2 > n_1$ such that $x_{e_{n_2}}^{n_2} \in B_2$. Thus we get soft subsequence $\left\{x_{e_{n_k}}^{n_k}\right\}$.

Now, let us show that the sequence $\left\{x_{e_{n_k}}^{n_k}\right\}$ is a Cauchy soft subsequence. For arbitrary $\tilde{\varepsilon} > \tilde{0}$, we choose number $\tilde{k}_0$ such that $\frac{\tilde{1}}{k_0} < \tilde{\varepsilon}$. Then since $x_{e_{n_k}}^{n_k}, x_{e_{n_m}}^{n_m} \in B_{k_0}$, for each $k, m \geq k_0$,

$$\tilde{d}\left(x_{e_{n_k}}^{n_k}, x_{e_{n_m}}^{n_m}\right) < \frac{\tilde{1}}{k_0} < \tilde{\varepsilon}$$

is satisfied. This means that $\left\{x_{e_{n_k}}^{n_k}\right\}$ is Cauchy soft subsequence.

$\Leftarrow$ Assume that $\left(\tilde{X},\tilde{d},E\right)$ is not totally bounded soft metric space. That means, $\left(\tilde{X},\tilde{d},E\right)$ has not finite $\tilde{\varepsilon}_0$ – net, for some $\tilde{\varepsilon}_0 > \tilde{0}$. Let $x_{e_1}^1 \in \tilde{X}$ be an arbitrary soft point. Then there can be found a soft point $x_{e_2}^2 \in \tilde{X}$ such that $\tilde{d}\left(x_{e_1}^1, x_{e_2}^2\right) \geq \tilde{\varepsilon}_0$. Since the soft set $\left\{x_{e_1}^1, x_{e_2}^2\right\}$ is not $\tilde{\varepsilon}_0$ – net, there is a soft point $x_{e_3}^3 \in \tilde{X}$ such that

$$\tilde{d}\left(x_{e_1}^1, x_{e_3}^3\right) \geq \tilde{\varepsilon}_0, \quad \tilde{d}\left(x_{e_2}^2, x_{e_3}^3\right) \geq \tilde{\varepsilon}_o.$$

Thus we constitute a soft sequence $\{x^k_{e_k}\}$ such that $\tilde{d}(x^i_{e_i}, x^j_{e_j}) \geq \tilde{\varepsilon}_0$, for all $i, j$. It is clear that $\{x^k_{e_k}\}$ has not a Cauchy soft subsequence which contradicts the fact that every soft sequence has a Cauchy soft subsequence given in assumption.

**Theorem 3.9.** Let $(\tilde{X}, \tilde{d}, E)$ be a soft metric space. $(\tilde{X}, \tilde{d}, E)$ is soft sequentially compact metric space if and only if $(\tilde{X}, \tilde{d}, E)$ is soft complete and soft totally bounded.

**Proof.** $\Rightarrow$ Let $(\tilde{X}, \tilde{d}, E)$ be a soft sequentially compact metric space. Then every soft sequence $\{x^n_{e_n}\}$ has a soft subsequence that converges in $\tilde{X}$. Since the soft subsequence is a Cauchy sequence, then by Theorem 3.8, $(\tilde{X}, \tilde{d}, E)$ is totally bounded. If $\{x^n_{e_n}\}$ is a Cauchy soft sequence in $(\tilde{X}, \tilde{d}, E)$ and $\{x^n_{e_n}\}$ has a convergent soft subsequence, then it is also convergent.

$\Leftarrow$ Let $(\tilde{X}, \tilde{d}, E)$ be a soft complete and totally bounded metric space and $\{x^n_{e_n}\}$ be an arbitrary soft sequence. Since $(\tilde{X}, \tilde{d}, E)$ is totally bounded, $\{x^n_{e_n}\}$ has a Cauchy soft subsequence. Since $(\tilde{X}, \tilde{d}, E)$ is soft complete, the Cauchy soft subsequence converges. Then $(\tilde{X}, \tilde{d}, E)$ is soft sequentially compact metric space.

**Definition 3.10.** Let $(\tilde{X}, \tilde{d}, E)$ be a soft metric space and a family $U$ be a soft open cover of the space $(\tilde{X}, \tilde{d}, E)$. A number $\tilde{\varepsilon} > \tilde{0}$ is called Lebesque number of $U$, if there exists $(F, E) \in U$ such that $B(x_e, \tilde{\varepsilon}) \subset (F, E)$, for all soft point $x_e \in \tilde{X}$.

**Proposition 3.11.** If $(\tilde{X}, \tilde{d}, E)$ is a soft sequentially compact metric space then every soft open cover from $\tilde{X}$ has a Lebesque number.

**Proof.** Assume that soft open cover $U$ has not a Lebesque number. Then for any $\tilde{n}$, there exists $x^n_{e_n}$, for each $(F, E) \in U$, where $B\left(x^n_{e_n}, \dfrac{\tilde{1}}{n}\right) \subset (F, E)$ is not satisfied. Thus we obtain a soft sequence $\{x_{e_n}^n\}$ satisfying above condition. Since $(\tilde{X},\tilde{d},E)$ is soft sequentially compact metric space, the soft sequence $\{x_{e_n}^n\}$ has a soft subsequence $\{x_{e_{n_k}}^{n_k}\}$ converging to $x_e$. Let $x_e \in (F,E) \in U$. Since $(F,E)$ is soft open set, there is soft open ball $B\left(x_e, \dfrac{\tilde{2}}{m}\right)$ such that $B\left(x_e, \dfrac{\tilde{2}}{m}\right) \subset (F,E)$. Also since the soft subsequence $\{x_{e_{n_k}}^{n_k}\}$ converges to $x_e$, there is a number $k_0$, such that $x_{e_{n_k}}^{n_k} \in B\left(x_e, \dfrac{\tilde{2}}{m}\right)$, whenever $k \geq k_0$. We take the number $k \geq k_0$ as $n_k \geq m$. Then

$$B\left(x_{e_{n_k}}^{n_k}, \dfrac{\tilde{1}}{n_k}\right) \subset B\left(x_e, \dfrac{\tilde{2}}{m}\right) \subset (F,E) \in U$$

is obtained and this contradicts with the choice of the soft point $x_{e_{n_k}}^{n_k}$.

**Theorem 3.12.** Let $(\tilde{X},\tilde{d},E)$ be a soft metric space. Then the following statements are equivalent:

**a)** $(\tilde{X},\tilde{d},E)$ is soft compact.

**b)** $(\tilde{X},\tilde{d},E)$ is sequentially compact.

**Proof.** $a \Rightarrow b$ Let $(\tilde{X},\tilde{d},E)$ be a soft compact metric space. Consider $(\tilde{X},\tilde{d},E)$ is not a soft sequentially compact metric space. Then there is an infinite soft set $\tilde{A}$ which does not have a cluster point in $(\tilde{X},\tilde{d},E)$. Thus, there is a soft number $r_{x_e}$ such that $B(x_e, r_{x_e}) \cap \tilde{A} = \{x_e\}$, for all $x_e \in \tilde{A}$. A family $\{B(x_e, r_{x_e})\}_{x_e \in \tilde{A}} \cup \tilde{A}^c$ is a soft open cover in $(\tilde{X},\tilde{d},E)$ and this soft open cover has not a finite soft subcover. But this implies $(\tilde{X},\tilde{d},E)$ is not soft compact, which is a contradiction.

$b \Rightarrow a$ Let $(\tilde{X},\tilde{d},E)$ be a soft sequentially compact and $U$ be any soft open cover in $(\tilde{X},\tilde{d},E)$. Then by Proposition 3.11, $U$ has a Lebesque number $\tilde{\varepsilon} > \tilde{0}$. Since $(\tilde{X},\tilde{d},E)$ is totally bounded, $(\tilde{X},\tilde{d},E)$ has finite $\dfrac{\tilde{\varepsilon}}{3}-$ net as $\{x_{e_1}^1, x_{e_2}^2, ..., x_{e_n}^n\}$. For each $k = 1,2,...,n$

$$\tilde{d}\left(B\left(x_{e_k}^k, \frac{\tilde{\varepsilon}}{3}\right)\right) \leq \frac{2\tilde{\varepsilon}}{3} < \tilde{\varepsilon}$$

is satisfied. Then we obtain a soft set $B\left(x_{e_k}^k, \frac{\tilde{\varepsilon}}{3}\right) \subset (F_k, E) \in U$. Since $\tilde{X} = \bigcup_{k=1}^n B\left(x_{e_k}^k, \frac{\tilde{\varepsilon}}{3}\right) \subset \bigcup_{k=1}^n (F_k, E)$, $\left(\tilde{X}, \tilde{d}, E\right)$ is a soft compact metric space.

**Definition 3.13.** Let $\left(\tilde{X}, \tilde{d}_1, E_1\right)$ and $\left(\tilde{Y}, \tilde{d}_2, E_2\right)$ be two soft metric spaces. The mapping $(f, \varphi): \left(\tilde{X}, \tilde{d}_1, E_1\right) \to \left(\tilde{Y}, \tilde{d}_2, E_2\right)$ is called soft uniformly continuous mapping if given any $\tilde{\varepsilon} > \tilde{0}$, there exists a $\tilde{\delta} > \tilde{0}$ ($\tilde{\delta}$ depending only on $\tilde{\varepsilon}$) such that for any soft point $x_e, y_{e'} \in \tilde{X}$ when $\tilde{d}_1(x_e, y_{e'}) < \tilde{\delta}$ satisfied then $\tilde{d}_2(f(x)_{\varphi(e)}, f(y)_{\varphi(e')}) < \tilde{\varepsilon}$.

**Proposition 3.14.** If $(f, \varphi): \left(\tilde{X}, \tilde{d}_1, E_1\right) \to \left(\tilde{Y}, \tilde{d}_2, E_2\right)$ is soft uniformly continuous mapping, then $f: (X, d_{1e}) \to (Y, d_{2\varphi(e)})$ is uniformly continuous mapping, for each $e \in E$.
**Proof.** The proof is straightforward.

The converse of the Proposition 3.14. above may not be true in general.

**Theorem 3.15.** If $(f, \varphi): \left(\tilde{X}, \tilde{d}_1, E_1\right) \to \left(\tilde{Y}, \tilde{d}_2, E_2\right)$ is soft continuous mapping and $\left(\tilde{X}, \tilde{d}_1, E_1\right)$ is soft sequentially compact metric space, then $(f, \varphi)$ is soft uniformly continuous mapping.
**Proof.** Since $\left(\tilde{X}, \tilde{d}_1, E_1\right)$ is soft sequentially compact metric space, it is also soft compact metric space. For any $\tilde{\varepsilon} > \tilde{0}$, since $(f, \varphi)$ is soft continuous, for any soft point $x_e$ there exists a soft number $\tilde{\delta}(x_e) > \tilde{0}$ such that for every soft point $y_{e'}$ that satisfies the condition $\tilde{d}_1(x_e, y_{e'}) < 2\tilde{\delta}(x_e)$ we have $\tilde{d}_2(f(x)_{\varphi(e)}, f(y)_{\varphi(e')}) < \frac{\tilde{\varepsilon}}{2}$. Then the family $U = \left\{B\left(x_e, \tilde{\delta}(x_e)\right)\right\}_{x_e \in \tilde{X}}$ is a soft open cover in $\left(\tilde{X}, \tilde{d}_1, E_1\right)$. Since $\left(\tilde{X}, \tilde{d}_1, E_1\right)$ is soft compact, this soft open cover has a finite soft subcover as $\left\{B\left(x_{e_1}^1, \tilde{\delta}(x_{e_1}^1)\right), \ldots, B\left(x_{e_n}^n, \tilde{\delta}(x_{e_n}^n)\right)\right\}$. We take

$$\tilde{\delta} = \min\left\{\tilde{\delta}(x_{e_1}^1),...,\tilde{\delta}(x_{e_n}^n)\right\}.$$

Now consider only two soft points $y_e, z_{e'} \in \tilde{X}$ such that $\tilde{d}_1(y_e, z_{e'}) < \tilde{\delta}$. Assume that $y_e \in B\left(x_{e_i}^i, \tilde{\delta}(x_{e_i}^i)\right)$, $1 \leq i \leq n$. Then $\tilde{d}_1(y_e, x_{e_i}^i) < \tilde{\delta}(x_{e_i}^i)$ and

$$\tilde{d}_1(z_{e'}, x_{e_i}^i) \leq \tilde{d}_1(z_{e'}, y_e) + \tilde{d}_1(y_e, x_{e_i}^i) < \tilde{\delta} + \tilde{\delta}(x_{e_i}^i) \leq 2\tilde{\delta}(x_{e_i}^i)$$

is satisfied. Since $(f, \varphi)$ is soft continuous at soft point $x_{e_i}^i$,

$$\tilde{d}_2\left(f(y)_{\varphi(e)}, f(x^i)_{\varphi(e_i)}\right) < \frac{\tilde{\varepsilon}}{2} \text{ and } \tilde{d}_2\left(f(z)_{\varphi(e')}, f(x^i)_{\varphi(e_i)}\right) < \frac{\tilde{\varepsilon}}{2}.$$

Thus

$$\tilde{d}_2\left(f(y)_{\varphi(e)}, f(z)_{\varphi(e')}\right) \leq \tilde{d}_2\left(f(y)_{\varphi(e)}, f(x^i)_{\varphi(e_i)}\right) + \tilde{d}_2\left(f(x^i)_{\varphi(e_i)}, f(z)_{\varphi(e')}\right) < \tilde{\varepsilon}.$$

As a result, $(f, \varphi)$ is a soft uniformly continuous mapping.

## 5. Conclusion

The soft set theory was proposed by Molodtsov offers a general mathematical tool for dealing with uncertain and vague objects. Many researchers have contributed towards the topologization of soft set theory. This study contributes some important theorems on soft sequential compact metric spaces. Later we give the concepts of $\tilde{\varepsilon}-$ net and totally bounded soft metric space. We continue to investigate some important properties of totally bounded soft metric spaces. Finally we introduce the concepts of Lebesque number and soft uniformly continuous mapping and investigate some theorems in detail.